\documentclass[11pt]{amsart}
\usepackage{amssymb}
\usepackage{xcolor,rotating,epic,eepic}
\usepackage{hhline}
\usepackage{fancyhdr}
\newcommand{\R}{\mathbb{R}}
\newcommand{\N}{\mathbb{N}}
\newcommand{\Z}{\mathbb{Z}}
\newcommand{\C}{\mathbb{C}}

\def\Log{\mathrm{Log} }

\def\endproof{\quad \hfill$\blacksquare$\vspace{0.15cm}\\}
\newcommand{\ds}{\displaystyle}

\newtheorem{prop}{ Proposition}

\newtheorem{cor}{Corollary}
\newtheorem{theorem}{Theorem}
\everymath{\displaystyle}

\begin{document}
\title[$q$-versions of Ramanujan master theorem]{On some $q$-versions of Ramanujan master theorem }

\author{Ahmed Fitouhi~~ \qquad \& \qquad~~Kamel Brahim\qquad \& \quad N\'eji Bettaibi}

\address{ Faculty of Sciences of Tunis. University of Tunis El Manar, 1060 Tunis, Tunisia}\email{
ahmed.fitouhi@fst.rnu.tn}
\address{
Faculty of Sciences of Tunis. University of Tunis El Manar, 1060 Tunis, Tunisia}\email{ Kamel.Brahim@ipeit.rnu.tn}
\address{ Department of Mathematics, College of Science,
Qassim University, Buraidah, Saudi Arabia}
\email{neji.bettaibi@ipein.rnu.tn} 

\maketitle
 \begin{abstract} In this paper, we state some $q$-analogues of the famous Ramanujan's Master Theorem. As applications, some values of   Jackson's $q$-integrals involving  $q$-special functions are computed.
 \end{abstract}
\section{Introduction}
In his First Quarterly Report to the Board of studies of the university of Madras in 1913, S. Ramanujan introduced a technique, known as Ramanujan Master Theorem (RMT), which provides an analytic expression for the Mellin transform of a function.\\
The theorem is asserted as follows: if a function $f$ can be expanded around $x=0$ in a power series in the form
$$f(x)= \sum_{n=0}^\infty (-1)^n \varphi (n)x^n$$
then
\begin{equation}\label{0.1}
\int_0^\infty x^{s-1}f(x)dx=\frac{\pi}{\sin(\pi s)}\varphi(-s).
\end{equation}
By replacing  $\varphi(s)$ by $\frac{\varphi(s)}{\Gamma(s+1)}$, we obtain the following equivalent version
\begin{equation}
\int_0^\infty x^{s-1}\left(  \sum_{n=0}^\infty (-1)^n \frac{\varphi (n)}{n!}x^n \right)dx=\Gamma(s)\varphi(-s).
\end{equation}
This theorem  was used by Ramanujan himself   as a tool in computing definite integrals and infinite series. In fact, most of the examples given by Ramanujan by applying his master theorem turn out to be correct, but this later can not hold without additional assumptions, as one can easily see from the example $\varphi(s)=\sin(\pi s)$. In \cite{GHH}, Hardy stated a rigorous reformulation of the  Ramanujan's Master Theorem. Using the residue theorem, Hardy proved the following alternative formulation:
\begin{theorem}\label{k1}
Let $f$ be a function having  around $x=0$ the expansion,
$$f(x)= \sum_{n=0}^\infty (-1)^n \varphi (n)x^n,$$
with $\varphi(z)$ is an analytic function defined on a half-plane :
$$H(\delta)=\{z\in \mathbb{C} : \Re(z)\geq-\delta \}$$
for some $0<\delta<1$.  Suppose that, for some $A<\pi$ and $P,\ C\in \R$, $\varphi$ satisfies the growth condition
$$|\varphi(v+iw)|<C e^{Pv+ A|w|},\;\; v+iw\in H(\delta).$$
Then for all  $0<Re s <\delta $, we have
\begin{equation}\label{ext}
\int_0^{\infty} x^{s-1}f(x)dx=\frac{\pi}{\sin (\pi x)} \varphi(-s).
\end{equation}
\end{theorem}
Moreover, using the Mellin's inversion formula, Hardy proved that under the conditions of the theorem, the function $f$ can be extended to $]0, \infty[$ and the formula (\ref{ext})  holds for the extension of $f$ on $]0, \infty[$. The condition $\delta<1$, ensures convergence of the integral in (\ref{ext}) and analytic continuation may be employed to valid  (\ref{ext}) to a large strip in which the integral converges.  \\

In view of the applications of the Ramnujan  Master Theorem in computing definite integrals and infinite series, and in the interpolation field,   Ismail and Stanton   in \cite{ID}  asked, in an open question,  about  the existence of  a $q$-analogue of this famous  theorem. Since the authors introduced and studied in \cite{FNK} a $q$-analogue of the Mellin transform, this question arouses their attention and the purpose of the present paper is to give an affirmative answer. We give two $q$-versions of the Ramnujan  Master Theorem. The proof of the  first $q$-version  is based on a  Mera's theorem \cite{Mera},  which gives a link between the Mellin transform and its $q$-analogue. For the proof of the second $q$-version, we use the  residue theorem by taking advantage of the asymptotic behavior of the $q$-Gamma function. Each of the two $q$-versions is followed by some consequences and examples.\\

In this paper, we adopt the following shot: Section 2 is devoted to present the elements of the quantum calculus making this paper independent. In Section 3,  we recall some properties of the $q$-Mellin transform   and we  give new ones.
Sections 4 and 5 are the main sections of the paper. They are devoted to state the two $q$-versions of the RMT and to give some of their  derivative results.  Finally, in Section 6, we apply  the results to compute some usually $q$-integrals.

\section{Elements of quantum calculus}
Throughout this text the fixed parameter of deformation is the real
$q$ selected such that $0<q<1$. To make this paper self containing  we recall the useful $q$-notions. We refer to the books \cite{GR} and \cite{Kac}, for the notations, definitions and
properties of the $q$-shifted factorials,  the basic hypergeometric functions and the Jackson's $q$-integrals.
For instance,

  {$\bullet$} For $a\in \C$, the $q$-shift factorials are defined by\\
$\ds (a;q)_0=1,~~ (a;q)_n={\prod_{k=0}^{n-1}(1-aq^k)}, ~~
n=1,2,\ldots, \infty .$\\

{$\bullet$} $\ds [x]_q=\frac{1-q^x}{1-q},~~x\in\C$~,~~~~$\ds
[n]_q!=\frac{(q;q)_n}{(1-q)^n},~~n\in\N.$\\

{$\bullet$} The
 $q$-hypergeometric series   is defined by\\
{\small $\ds _r\varphi_s(a_1   ...  a_r;b_1  ... b_s;q,z)=  ~~
\sum_{k=0}^\infty \frac{(a_1,a_2,...,a_r;q)_k}{(b_1,b_2,...,b_s;q)_k(q;q)_k}\left((-1)^kq^{\frac{k(k-1)}{2}}\right)^{1+s-r}z^k,$\\
where  $r,s\in \N$, $a_1,..., ~~a_r,~b_1,...,~~b_s\in \C$ such that
$b_1,...,~b_s\neq 1,~q^{-1},~q^{-2},...,$} and  \\
$(a_1,a_2,...,a_r;q)_k=(a_1;q)_k(a_2;q)_k...(a_r;q)_k.$\\

$\bullet$  The Jackson's $q$-integrals from $0$ to $a$ and from $0$
to $\infty$ are given by
{\small \begin{eqnarray*} \int_0^{a}{f(x)d_qx}
=(1-q)a\sum_{n=0}^{\infty}{f(aq^n)q^n},
~~~~\int_0^{\infty}{f(x)d_qx}
=(1-q)\sum_{n=-\infty}^{\infty}{f(q^n)q^n}.
\end{eqnarray*}

{$\bullet$} The improper $q$-integral is defined in the following way\\
$$\int_0^{\infty\over A} f(x)d_qx= (1-q)\sum_{n=-\infty}^\infty f\left(\frac{q^n}{A}\right)\frac{q^n}{A}.$$

{$\bullet$} Two $q$-analogues of the  exponential function are given by (see \cite{Kour})
\begin{equation}
E_q^z = \sum_0^\infty
q^{n(n-1)}\frac{z^n}{[n]_n!}=(-(1-q)z;q)_{\infty}
\end{equation}
and
\begin{equation}
e_q^z=\sum_0^\infty \frac{z^n}{[n]_n!}
={1\over{((1-q)z;q)_{\infty}}},~~~~~~~~|z|<\frac{1}{1-q}.
\end{equation}

{$\bullet$} The $q$-gamma function is defined by

\begin{equation}
 \Gamma_q (s)
={(q;q)_{\infty}\over{(q^s;q)_{\infty}}}(1-q)^{1-s} ,
\qquad s\not={0},{-1},{-2},\ldots,\\
\end{equation}
and has the following  $q$-integral
representations:
\begin{equation}
\Gamma_q(s)=\int_0^{\frac{\infty}{1-q}}t^{s-1}E_q^{-qt}d_qt= K_q(s)\int_0^{\frac{\infty}{1-q}}t^{s-1}e_q^{-t}d_qt,\\
\end{equation}
where
\begin{equation}\label{Kq}
\ds K_q(s)= \frac{(-q,-1,q)_{\infty}}{(-q^s, -q^{1-s};q)_{\infty}}.\\
\end{equation}

{$\bullet$} The third Jackson's $q$-Bessel function,  called also in some literature  Heine Exton Bessel function, is  defined by
$$ J_{\alpha}(z;q^2)
=\frac{z^{\alpha}}{(1-q^2)^{\alpha}\Gamma_{q^2}(\alpha+1)}
        {_1\varphi_1(0;q^{2\alpha + 2};q^2,q^2z^2)}.
$$

In the end of this section, we enunciate that we will use the following sets:

{$\bullet$}  $\mathbb{R}_{q}= \{\pm q^n: n\in \Z\},\qquad \qquad~~~\mathbb{R}_{q,+} = \{q^n: n\in \Z\}$.\\

{$\bullet$} For $a,
 b\in \mathbb{R},~~ (a<b)$, we write\\
\indent \qquad\qquad { $ \langle a; b \rangle = \left\{ s\in\C
~~:~~a<\Re(s)<b\right\}$.}

\section{$q$-Mellin transform}
Hjalmar Mellin introduced the Mellin transform for a suitable function $f$ defined over $]0,\infty[$ by
\begin{eqnarray*}
M(f)(s)=\int_0^{\infty} f(x)x^{s-1}dx.
\end{eqnarray*}
This integral is defined in some open strip {\bf$\langle \alpha_{f};\beta_{f} \rangle$} called fundamental band depending on the behavior of $f$ at $0$ and $\infty$ .\\
In \cite{FNK}, Fitouhi et al.  introduced and studied  a $q$-analogue of the Mellin transform, called  $q$-Mellin transform, by the use of the Jackson's $q$-integral:  for a  function  $f$ defined on  $\mathbb{R}_{q,+}$,
\begin{eqnarray}
 M_q(f)(s)=M_q[f(t)](s)=\int_0^{\infty}t^{s-1}f(t)d_qt=(1-q)\sum_{n=-\infty}^\infty f(q^n)q^{ns}.
\end{eqnarray}

The function  $M_q(f)$ tends to $M(f)$ when $q$ tends to $1^-$, furthermore it is analytic on a band {\bf$\langle \alpha_{q,f};\beta_{q,f} \rangle$}, called the fundamental strip of $M_q(f)$,   and it has the important property namely  $\ds \frac{2i\pi}{\mathrm{Log} (q)}$-periodic.\\

Other important properties and applications are proved rigorously in Fitouhi et al \cite{FNK},  among which we just cite

\begin{eqnarray}
M_q[f\left(\frac{1}{t}\right)](s)=M_q(f)(-s),\end{eqnarray}
\begin{eqnarray}
M_q[f\left(x^\rho\right)](s)=\left[ \frac{1}{\rho}\right]_{q^\rho}M_{q^\rho}(f)\left(\frac{s}{\rho}\right),\quad \rho>0\end{eqnarray}
and
\begin{eqnarray}
M_q[t^xf\left(t\right)](s)=M_q(f)(s+x).\end{eqnarray}

In \cite{FNK}, the $q$-Mellin inversion formula is established and proved. We shall give, here,  an another  proof based on  the periodicity of this $q$-Mellin transform    and   on the  Fourier series.

\begin{theorem} Let $f$  be a function defined on $\mathbb{R}_{q,+}$ and $M_q(f)$ be its $q$-Mellin transform with {\bf$\langle \alpha_{q,f}; \beta_{q,f} \rangle$} as fundamental band. Then, for   $c\in
\langle \alpha_{q,f}; \beta_{q,f}\rangle\cap \mathbb{R},$ we have\\
\begin{eqnarray}
\forall x \in \mathbb{R}_{q,+}, \quad f(x) =
\frac{\mathrm{Log}(q)}{2i\pi (1-q)}
\int_{c-\frac{i\pi}{\mathrm{Log}(q)}}^{c+\frac{i\pi}{\mathrm{Log}(q)}}
M_q(f)(s)x^{-s}ds.
\end{eqnarray}
\end{theorem}

\proof ~~\\ Let  $c$ be a real in $\langle\alpha_{q,f};\beta_{q,f}\rangle$. Put
$$g(t)=M_q(f)(c+it)=(1-q)\sum_{n=-\infty}^{\infty}f(q^n)q^{nc}e^{{int}\Log
q}.$$
Since   $M_q(f)$ is analytic on the band {\bf$\langle\alpha_{q,f}; \beta_{q,f} \rangle$} and $\frac{2i\pi}{Log q}$-periodic, then $g$ is continuous and $\frac{2\pi}{\Log (q)}$-periodic, on $\R$.
Hence, $g$ can be expanded in Fourier series and the convergence of the series is uniform. That is
$$ g(t)=\sum_{n=-\infty}^\infty C_n(g) e^{int \Log q},$$ where the coefficients $C_n(g)$ are given by

$$C_n(g)=\frac{\Log q}{2\pi}\int_{-\frac{\pi}{\Log q}}^{\frac{\pi}{\Log
q}}M_q(f)(c+it)e^{-int \Log q} dt.$$
From the uniqueness of the development of the Fourier series, we get for all $n\in \Z$,
 \begin{eqnarray*} C_n(g)=(1-q)f(q^n)q^{nc}&=&\frac{\Log q}{2\pi}\int_{-\frac{\pi}{\Log q}}^{\frac{\pi}{\Log
q}}M_q(f)(c+it)e^{-int \Log q} dt\\ &=& \frac{\Log q}{2\pi}\int_{-\frac{\pi}{\Log q}}^{\frac{\pi}{\Log
q}}M_q(f)(c+it)q^{-int} dt\\ &=& \frac{\Log q}{2i\pi}\int_{c-\frac{i\pi}{\Log q}}^{c+\frac{i\pi}{\Log
q}}M_q(f)(s)q^{n(c-s) } ds.\end{eqnarray*}
Thus, for all  $x=q^n\in \mathbb{R}_{q,+}$, we have
 $$ f(x)=\frac{\Log q}{2i\pi(1-q)}\int_{c-\frac{i\pi}{\Log q}}^{c+\frac{i\pi}{\Log
q}}M_q(f)(s)x^{-s} ds.$$
\endproof
The following result characterizes the functions $F$ which are the $q$-Mellin transforms  of functions defined on $\mathbb{R}_{q,+}$.
\begin{theorem}\label{car}
Let $a$ and $b$ be two real numbers  such that  $a<b$ and $F$ be an analytic and $\frac{2i\pi}{\Log q}$-periodic function  on the
band $\langle a;b\rangle$. Then,  $F $ is the $q$-Mellin transform of  the function
$f$ defined on $\mathbb{R_{q,+}}$  by
$$ f(x)=\frac{\Log q}{2i\pi(1-q)}\int_{c-\frac{i\pi}{\Log q}}^{c+\frac{i\pi}{\Log
q}}F(s)x^{-s} ds,$$
where $c$ is an arbitrary real in  $\langle a;b\rangle$.
\end{theorem}

\proof ~~
For $a<c<b$, we consider the function $h$ defined on $\R$ by
$$h(t)=F(c+it).$$
Evidently,  the function $h$ is continuous, differentiable and $\frac{2\pi}{\Log q}-$periodic on $\mathbb{R}$, so  it can be expanded in Fourier series as
$$h(t)=\sum_{-\infty}^\infty C_n(h)e^{int \Log q},$$
where
\begin{eqnarray*}C_n(h)&=&\frac{\Log q}{2\pi}\int_{-\frac{\pi}{\Log q}}^{\frac{\pi}{\Log
q}}h(t) e^{-in t \Log q}dt= \frac{\Log q}{2\pi}\int_{-\frac{\pi}{\Log q}}^{\frac{\pi}{\Log
q}}F(c+it) q^{-in t}dt\\
&=& \frac{\Log q}{2i\pi}\int_{c-\frac{i\pi}{\Log q}}^{c+\frac{i\pi}{\Log
q}}F(s) q^{-n(s-c)}ds=q^{nc}\frac{\Log q}{2i\pi}\int_{c-\frac{i\pi}{\Log q}}^{c+\frac{i\pi}{\Log
q}}F(s) q^{-ns}ds\\&=& q^{nc}(1-q)f(q^n),
\end{eqnarray*}
where the function $f$ is defined on $\R_{q,+}$ by
 $$f(q^n)=\frac{\Log q}{2i\pi(1-q)}\int_{c-\frac{i\pi}{\Log q}}^{c+\frac{i\pi}{\Log
q}}F(s) q^{-ns}ds,  n\in \mathbb{Z}.$$
Then, for $s=c+it$, we have
 \begin{eqnarray*}F(s)&=&F(c+it)= h(t)=(1-q)\sum_{n=-\infty}^\infty f(q^n)q^{nc} e^{int \Log q}\\&=&(1-q)\sum_{n=-\infty}^\infty f(q^n)q^{n(c+it)} =(1-q)\sum_{n=-\infty}^\infty f(q^n)q^{n(s)}\\&=& \int_{0}^\infty x^{s-1}f(x)d_qx = M_q(f)(s).\end{eqnarray*}
Finally, we remark that since $F$ is $\frac{2i\pi}{\Log q}-$ periodic and analytic on $\langle a;b\rangle$, then thanks to the Cauchy's theorem, the function $f$ does not depend on the choice of the real  $c$ in $\langle a;b\rangle$.\endproof

\section{First $q$-Version of Ramanujan master theorem}

To establish the main   result  of this section, we need the following result, which is recently stated by  M. Mera theorem in \cite{Mera} and gives a  relationship between  the classical Mellin and the   $q$-Mellin transforms.

\begin{theorem}
Let $f$ a function defined on $(0,+\infty)$ and $M(f)$ be  its
Mellin transform with  $\langle i_0,j_0\rangle$ as fundamental
strip. For $s\in \langle i_0,j_0\rangle$ and  $N\in \mathbb{N}\cup
\{0\}$,    we put
\begin{equation}
M_q^N(f)(s)=-\frac{1-q}{ \Log q}\sum_{|n|\leq N}M(f)\left(s+i\frac{2\pi n}{ \Log q}\right).
\end{equation}
Suppose that for a fixed $\sigma$, such that  $-i_0<\sigma<-j_0$,  the sum $M_q^N(f)(\sigma +it)$ converges absolutely and uniformly  for $t$ in any compact subset of $\mathbb{R}$ as $N\rightarrow \infty$.
Then
$$\lim_{N \rightarrow \infty}M_q^N(f)(\sigma + it)=M_q(f)(\sigma + it); t\in \mathbb{R}.$$\\
\end{theorem}

\noindent{\bf Example}\\ Consider the function $f$ defined on $(0,+\infty)$ by
$$f(x)= \frac{1}{1+x}.$$
We have $$ M(f)(s)= \Gamma(s)\Gamma(1-s)=\frac{\pi}{\sin (\pi s)}, \; \; 0<\Re(s)<1 $$ and for all $N\in\N\cup \{0\}$,
\begin{equation}\label{eq10}
M_q^N(f)(s)= -\frac{1-q}{\Log q}\sum_{|n|\leq N}\frac{\pi}{\sin \left(\pi \left(s+\frac{2i\pi n}{\Log q}\right)\right)}.
\end{equation}
Using the relation $\Gamma(s)= \frac{\Gamma(s+2)}{s(s+1)}$  one can deduce that for $s=\sigma +it\in \langle 0;1\rangle $, we have
$$\left|\frac{\pi}{\sin \left(\pi(s+\frac{2i\pi n}{\Log q})\right)}\right|\leq \frac{\Gamma(\sigma+2)\Gamma(1-\sigma)}
{\sigma^2+(t+\frac{2\pi n}{\Log q})^2}.$$
So, the sum $M_q^N(f)(\sigma +it)$ converges absolutely and uniformly for $t$ in any compact subset of $\mathbb{R}$. Then, from the previous theorem, we have   for $0< \Re (s)<1$,

$$\lim_{N\rightarrow \infty}M_q^N\left[\frac{1}{1+x}\right](s)=M_q\left[\frac{1}{1+x}\right](s)=\frac{\Gamma_q(s)\Gamma_q(1-s)}{K_q(s)}.$$\\

The following theorem is the main result of this section. It is the object of the first $q$-version of the Ramanujan master theorem. The proof is based on the previous Mera's theorem.
\begin{theorem}\label{Mq1}
 Let $\varphi$ be a function satisfying  the conditions of Theorem 1 \label{th1} and   $\frac{2i\pi}{\Log q}-$ periodic. Then for $0<\Re(s)<\delta $, we have
 \begin{equation}\label{ext1}
\int_0^\infty x^{s-1}\left(\sum_{k=0}^\infty \varphi(k) (-x)^k\right)d_qx = \frac{\Gamma_q(s)\Gamma_q(1-s)}{K_q(s)}\varphi(-s),
\end{equation}
where $K_q(.)$ is  defined by (\ref{Kq}).
\end{theorem}
\proof ~~\\
Put
$$f(x)= \sum_{k=0}^\infty (-1)^k \varphi(k) x^k.$$
Using the Ramanujan master theorem, we get  \\
$$ M(f)(s)=\varphi(-s)\frac{\pi}{\sin(\pi s) }, \; \; 0<\Re(s)<\delta.$$
Then, thanks to the   Mera's result,  the periodicity of $\varphi$ and the relation (\ref{eq10}), we have for all $N\in \N$,
\begin{eqnarray*}M_q^N(f)(s)&=&  -\frac{1-q}{\Log q}\sum_{|n|\leq N} M(f)\left(s+\frac{2i\pi n}{\Log q}\right)\\ &=&  -\frac{1-q}{\Log q}\sum_{|n|\leq N}\frac{\pi}{\sin\left( \pi\left(s+\frac{2i\pi n}{\Log q}\right)\right)} \varphi(-s)\\ &=& \varphi(-s) M_q^N\left[\frac{1}{1+x}\right](s). \end{eqnarray*}
Since for $0<\sigma<1$, the sum $M_q^N\left[\frac{1}{1+x}\right](\sigma +it)$ converges absolutely  and uniformly for $t$ in any compact of $\R$, then for $0<\sigma<\delta$, the sum $M_q^N(f)(\sigma +it)$ converges absolutely  and uniformly for $t$ in any compact of $\R$ and for $s=\sigma+it$, we have
\begin{eqnarray*}M_q(f)(s)&=& \lim_{N \rightarrow \infty}M_q^N(f)(s)=\varphi(-s) \lim_{N \rightarrow \infty}M_q^N\left[\frac{1}{1+x}\right](s)\\  &=&\varphi(-s)M_q\left[\frac{1}{1+x}\right](s)= \frac{\Gamma_q(s)\Gamma_q(1-s)}{K_q(s)}\varphi(-s).\end{eqnarray*}\endproof

\noindent{\bf Remark}\\Analytic continuation may be employed to valid  (\ref{ext1}) to a large strip in which the $q$-integral converges. So, in the remainder, we will not precise the strip of convergence. \\

\noindent{\bf Example}\\
The function $f:x\mapsto\frac{1}{(-x; q)_\infty}$ is defined on
$\Re(x)\geq 0$ and for $|x|<1$, we have
$$f(x)=\sum_{n=0}^\infty (-1)^n \varphi(n) x^n,$$
where   $\varphi(s)=\frac{(q^{s+1}; q)_\infty}{(q ; q)_\infty}$.\\
It is easy to see that the function $\varphi$ is entire, $\frac{2i\pi}{\Log q}-$periodic on $\C$ and for $\Re(s)>-\delta >-1$, we have
$$|\varphi(s)|\leq \frac{(-1 ; q)_{\infty}}{(q ; q)_{\infty}}.$$
So, for $0<\Re(s)<\delta$, we have
$$M_q(f)(s)=\frac{(q^{1-s} ; q)_{\infty}}{(q ; q)_{\infty}}\frac{\Gamma_q(s)\Gamma_q(1-s)}{K_q(s)}=(1-q)^s\frac{\Gamma_q(s)}{K_q(s)}.$$
Since the $q$-integral $M_q(f)$ converges for $\Re(s)>0$, then this relation can be extended to $\Re(s)>0$.
As  direct consequence, we get
$$\int_0^{\frac{\infty}{1-q}} \frac{x^{s-1}}{(-(1-q)x ; q)_{\infty}} d_q x=\frac{\Gamma_q(s)}{K_q(s)}.$$

\begin{prop}
Let $\varphi$ be a function satisfying the conditions of Theorem \ref{Mq1}, and put
$$f(x)=\sum_{n=0}^{\infty} (-1)^n\varphi(n) \frac{x^n}{(q; q)_n}.$$
Then
\begin{equation}
M_q(f)(s)= (1-q)^s\varphi(-s)\frac{\Gamma_q(s)}{K_q(s)}.
\end{equation}
\end{prop}
\proof ~~\\
Remark that
$$\forall n\in\N, \; \; \frac{\varphi(n)}{(q; q)_n}=\varphi(n) \frac{(q^{n+1}; q)_\infty}{(q; q)_\infty}$$
and put $$\psi(s)=\varphi(s) \frac{(q^{s+1}; q)_\infty}{(q; q)_\infty}.$$
If $\varphi$    satisfies the conditions of Theorem \ref{Mq1}, then the function $\psi$ is $\frac{2i\pi}{\Log( q)}-$ periodic and analytic on $H(\delta)$,
and for $s=v+iw\in H(\delta)$, we have $v\geq -\delta\geq -1$ and
$$|\psi(s)|=\left|\varphi(s) \frac{(q^{s+1}; q)_\infty}{(q; q)_\infty}\right|\leq C\frac{(-1; q)_\infty }{(q; q)_\infty}e^{Pv+A|w|}. $$
Then, $\psi$ satisfies the conditions of Theorem \ref{Mq1}, so  for $0<\Re(s)<\delta$, we have
\begin{eqnarray*}M_q(f)(s)=M_q\left[\sum_{n=0}^{\infty} (-1)^n\psi(n) x^n \right](s)=\frac{\Gamma_q(s)\Gamma_q(1-s)}{K_q(s)}\psi(-s)= (1-q)^s\varphi(-s)\frac{\Gamma_q(s)}{K_q(s)}. \end{eqnarray*}
\endproof
\begin{cor}\label{cor1}
Let $\varphi$ be an analytic function on $H(\delta)$ such that the function $s\mapsto (1-q)^s\varphi(s)$ satisfies the conditions of Theorem \ref{Mq1}. Then,
for $0<\Re(s)<\delta$, we have
\begin{equation}
M_q\left[\sum_{n=0}^{\infty} (-1)^n\varphi(n) \frac{x^n}{[n]_q!}\right]= \varphi(-s)\frac{\Gamma_q(s)}{K_q(s)}.
\end{equation}
\end{cor}

In the particular case,  $\frac{\Log (1-q)}{\Log q} \in \mathbb{Z}$, we have:
\begin{cor} If $\varphi$ is a function satisfying the conditions of Theorem \ref{Mq1}, then
for $0<\Re(s)<\delta$, we have
\begin{equation}
M_q\left[\sum_{n=0}^{\infty} (-1)^n\varphi(n) \frac{x^n}{[n]_q!}\right]= \varphi(-s)\frac{\Gamma_q(s)}{K_q(s)}.
\end{equation}In particular  we have
$$M_q\left[\frac{1}{(-(1-q)x ; q)_\infty}\right](s)=\frac{\Gamma_q(s)}{K_q(s)}.$$
\end{cor}

\begin{prop}\label{propal}
Let $\alpha>0$, $p\in\N\cup\{0\}$ and  $\varphi$ be an analytic function on $H(\delta)$ such that the function $s\mapsto (1-q^{\alpha})^s \varphi(s)$ is $\frac{2i\pi}{\Log\left(q^{\alpha}\right)}-$ periodic and satisfies the conditions of Theorem \ref{k1}.
Put
$$ f(x)= \sum_{n=0}^{\infty}(-1)^n\varphi(n)\frac{x^{\alpha n}}{[n+p]_{q^{\alpha}}!}.$$
Then,
\begin{equation}
M_q(f)(s)= \left[\frac{1}{\alpha}\right]_{q^{\alpha}} \varphi\left(- \frac{s}{\alpha}\right) \frac{\Gamma_{q^{\alpha}}\left(1-\frac{s}{\alpha}\right)}{\Gamma_{q^{\alpha}}\left(1-\frac{s}{\alpha}+p\right)}\frac{\Gamma_{q^{\alpha}}(\frac{s}{\alpha})}{K_{q^{\alpha}}(\frac{s}{\alpha})}
\end{equation}
\end{prop}
\proof ~~\\
On the one hand, we have
$$ f(x)= \sum_{n=0}^{\infty}(-1)^n\varphi(n)\frac{x^{\alpha n}}{[n+p]_{q^{\alpha}}!}= h(x^{\alpha}),$$
where
$$ h(x)= \sum_{n=0}^{\infty}(-1)^n\varphi(n)\frac{x^{n}}{[n+p]_{q^{\alpha}!}}.$$
So, using the properties of the $q$-Mellin transform, we get
\begin{equation}\label{Mk2}
 M_q(f)(s)= \left[\frac{1}{\alpha}\right]_{q^{\alpha}}M_{q^{\alpha}}(h)\left(\frac{s}{\alpha}\right).
 \end{equation}
On the other hand, since the function $ s\mapsto \frac{\varphi(s)}{\Gamma_{q^{\alpha}}(1+s+p)}=\frac{\left(q^{\alpha(s+p+1)};q^\alpha \right)_\infty}{\left(q^{\alpha};q^\alpha \right)_\infty}\left(1-q^\alpha \right)^{s+p}\varphi (s)$ satisfies the conditions of Theorem \ref{k1} and it is $\frac{2i\pi}{Log(q^{\alpha})}-$ periodic, then  by applying Theorem \ref{Mq1}, we obtain
$$ M_{q^{\alpha}}(h)(s)= \frac{\varphi(-s)}{\Gamma_{q^{\alpha}}(1-s+p)}\frac{\Gamma_{q^{\alpha}}(1-s)\Gamma_{q^{\alpha}}(s)}{K_{q^{\alpha}}(s)}.$$
Hence, the result follows  from the relation (\ref{Mk2}).\endproof

For the particular case $p=0$, we obtain the following result:
\begin{cor}
Let $\alpha>0$ and $\varphi$ be a function satisfying the condition of Theorem \ref{k1} such that the function $s\mapsto (1-q^{\alpha})^s \varphi(s)$ is $\frac{2i\pi}{\Log(q^{\alpha})}-$ periodic.
Then,
\begin{equation}
M_q\left[ \sum_{n=0}^{\infty}(-1)^n\varphi(n)\frac{x^{\alpha n}}{[n]_{q^{\alpha}}!}  \right](s)= \left[\frac{1}{\alpha}\right]_{q^{\alpha}} \varphi\left(- \frac{s}{\alpha}\right) \frac{\Gamma_{q^{\alpha}}\left(\frac{s}{\alpha}\right)}{K_{q^{\alpha}}\left(\frac{s}{\alpha}\right)}.
\end{equation}
\end{cor}
\begin{prop}
Let $\alpha\in\N$, $p\in\N\cup\{0\}$ and $\varphi$ be a function satisfying the condition of Theorem \ref{k1} and such that the function $s\mapsto (1-q)^{\alpha s} \varphi(s)$ is $\frac{2i\pi}{Log(q^{\alpha})}-$ periodic. Then, the $q$-Mellin transform of the function
$$ f(x)= \sum_{n=0}^{\infty}(-1)^n\varphi(n)\frac{x^{\alpha n}}{[\alpha n+p]_q!}$$
is given by
\begin{equation}
M_q(f)(s)= \left[\frac{1}{\alpha}\right]_{q^{\alpha}} \varphi\left(- \frac{s}{\alpha}\right) \frac{\Gamma_{q^{\alpha}}\left(1-\frac{s}{\alpha}\right)}{\Gamma_{q}\left(1-s+p\right)}\frac{\Gamma_{q^{\alpha}}\left(\frac{s}{\alpha}\right)}{K_{q^{\alpha}}\left(\frac{s}{\alpha}\right)}.
\end{equation}
\end{prop}
\proof ~~\\
Following the same steps of the proof of Proposition \ref{propal}  by replacing the function $ s\mapsto \frac{\varphi(s)}{\Gamma_{q}(1+ s+p)}$ by the function$ s\mapsto \frac{\varphi(s)}{\Gamma_{q}(1+\alpha s+p)}$, we get the desired result.\endproof

\section{Second $q$-version of the Ramanujan master theorem}
In this section, we use the the residue theorem, to prove a new $q$-version of the Ramanujan master Theorem.
\begin{theorem}\label{RMq2}
Let $\varphi$ be an analytic function defined on a half-plane
$$H(\delta)= \left\{z\in\C  :  \Re(z)\geq -\delta  \right\}$$
for some $0<\delta<1$. Suppose that $\varphi$ is $\frac{2i\pi}{Log(q)}$-periodic, and  there exist a continuous function $\varphi_1$ on $\R$ and real $r\in\R$ such that
\begin{equation}\label{maj} \forall z=u+iw\in H(\delta),\ \ |\varphi(u+iw)|\leq e^{r|u|}\varphi_1(w).\end{equation}
Then, for all $s\in \langle 0,\delta\rangle$, we have
\begin{equation}
\int_0^\infty x^{s-1}\left(\sum_{n=0}^\infty (-1)^nq^{\frac{n(n+1)}{2}}\varphi(n)x^n\right)d_qx= \Gamma_q(s)\Gamma_q(1-s)\varphi(-s).
\end{equation}
\end{theorem}
\proof~~\\
Let $\varphi$ be a function satisfying the conditions of the theorem 1 and $x=q^m\in\R_{q,+}$. Then the function $\psi$ defined by
$$\psi(s)=  \Gamma_q(s)\Gamma_q(1-s)\varphi(-s)x^{-s} $$
is $\frac{2i\pi}{Log(q)}$-periodic and analytic on $\langle 0,\delta\rangle$ and meromorphic on the half-plane $\Re(s)\leq \delta$, with simple poles at $0,-1,-2,-3,\dots$.
Furthermore, for all $n\in\N\cap \{0\}$, we have
\begin{eqnarray*}
Res(\psi, -n)&=& \lim_{s\rightarrow -n}(s+n)\psi(s)\\&=& \lim_{s\rightarrow -n}\left[(s+n)\Gamma_q(s)\Gamma_q(1-s)\varphi(-s)x^{-s}\right]\\
&=&\frac{1-q}{\Log(q)} q^{\frac{n(n+1)}{2}}(-1)^{n+1}\varphi(n)x^n.
\end{eqnarray*}
The conditions on the function $\varphi$ and the presence of the factor $q^{\frac{n(n+1)}{2}}$ show that the series $\sum_{n\geq 0}\frac{1-q}{\Log(q)} q^{\frac{n(n+1)}{2}}(-1)^{n+1}\varphi(n)x^n$ converges.\\
On the other hand, for $N\in\N$ and $0<\sigma<\delta$, the application of the residue theorem for the function $\psi$ using  the rectangular contour $C_N$  defined by the segments $\left[\sigma -\frac{i\pi}{\Log(q)}, \sigma +\frac{i\pi}{\Log(q)} \right]$, $\left[\sigma +\frac{i\pi}{\Log(q)}, -N+\frac{1}{2} +\frac{i\pi}{\Log(q)} \right]$, $\left[-N-\frac{1}{2} +\frac{i\pi}{\Log(q)}, -N-\frac{1}{2} -\frac{i\pi}{\Log(q)} \right]$\\  and $\left[-N-\frac{1}{2} -\frac{i\pi}{\Log(q)}, \sigma-\frac{i\pi}{\Log(q)} \right]$
yields
\begin{eqnarray}\label{Res}
-2i\pi \sum_{n=0}^N Res(\psi, -n)= \int_{C_N}\psi(s)ds,
\end{eqnarray}
where
\begin{eqnarray*}
 \int_{C_N}\psi(s)ds
&=& \int_{\sigma -\frac{i\pi}{\Log(q)}}^{ \sigma +\frac{i\pi}{\Log(q)}}\psi(s)ds+ \int_{\sigma +\frac{i\pi}{\Log(q)}}^{ -N-\frac{1}{2} +\frac{i\pi}{\Log(q)} }\psi(s)ds + \int_{-N-\frac{1}{2} +\frac{i\pi}{\Log(q)}}^{-N-\frac{1}{2} -\frac{i\pi}{\Log(q)} }\psi(s)ds\\  &+&\int_{-N-\frac{1}{2} -\frac{i\pi}{\Log(q)}}^{ \sigma-\frac{i\pi}{\Log(q)} }\psi(s)ds.
\end{eqnarray*}
The periodicity and the analyticity of the function $\psi$ imply that the sum of the second and the
fourth integrals is zero. Moreover, the third integral verifies
\begin{eqnarray*}
\int_{-N-\frac{1}{2} +\frac{i\pi}{\Log(q)}}^{-N-\frac{1}{2} -\frac{i\pi}{\Log(q)} }\psi(s)ds&=& i \int_{\frac{\pi}{\Log(q)}}^{-\frac{\pi}{\Log(q)} }\psi(-N-\frac{1}{2}+it)dt\\ &=&  i \int_{\frac{\pi}{\Log(q)}}^{-\frac{\pi}{\Log(q)} }\Gamma_q(-N-\frac{1}{2}+it)\Gamma_q( N+\frac{3}{2}-it) \varphi(N+\frac{1}{2}-it)x^{N+\frac{1}{2}-it} dt
\end{eqnarray*}
and
\begin{eqnarray*}
\left|\int_{-N-\frac{1}{2} +\frac{i\pi}{\Log(q)}}^{-N-\frac{1}{2} -\frac{i\pi}{\Log(q)} }\psi(s)ds\right|&\leq&
\int_{\frac{\pi}{\Log(q)}}^{-\frac{\pi}{\Log(q)} }\left|\Gamma_q(-N-\frac{1}{2}+it)\Gamma_q( N+\frac{3}{2}-it)\right| |\varphi(N+\frac{1}{2}-it)|x^{N+\frac{1}{2}} dt.
\end{eqnarray*}
But for all $t\in \left[\frac{\pi}{\Log(q)}, - \frac{\pi}{\Log(q)}\right] $  and  all $N\in \N$, we have
\begin{eqnarray*}
\Gamma_q(-N-\frac{1}{2}+it)\Gamma_q( N+\frac{3}{2}-it)&=& \frac{(1-q)^N}{(q^{-N-\frac{1}{2}+it};q)_N}\Gamma_q(-\frac{1}{2}+it) \frac{(q^{\frac{3}{2}-it};q)_N}{(1-q)^N}\Gamma_q( \frac{3}{2}-it)\\&=& \frac{(q^{\frac{3}{2}-it};q)_N}{(q^{-N-\frac{1}{2}+it};q)_N}\Gamma_q(-\frac{1}{2}+it) \Gamma_q( \frac{3}{2}-it).
\end{eqnarray*}
Since
\begin{eqnarray*}
|(q^{\frac{3}{2}-it};q)_N|=\prod_{k=0}^{N-1}|1-q^{\frac{3}{2}-it+k}|\leq \prod_{k=0}^{N-1}(1+q^{\frac{3}{2}+k})\leq (-q^{\frac{3}{2}};q)_\infty
\end{eqnarray*}
and
\begin{eqnarray*}
|(q^{-N-\frac{1}{2}+it};q)_N|&=&\prod_{k=0}^{N-1}|1-q^{-N-\frac{1}{2}+it+k}|\geq \prod_{k=0}^{N-1}(q^{-\frac{1}{2}-N+k}-1)=\prod_{k=0}^{N-1}q^{-\frac{1}{2}-N+k}(1-q^{\frac{1}{2}+N-k})\\ &\geq& q^{-\frac{N}{2}-\frac{N(N+1)}{2}}\frac{(q^{\frac{1}{2}};q)_{N+1}}{1-q^{\frac{1}{2}}}\geq q^{-\frac{N(N+2)}{2}}\frac{(q^{\frac{1}{2}};q)_{\infty}}{1-q^{\frac{1}{2}}},
\end{eqnarray*}
then
$$\left|  \Gamma_q(-N-\frac{1}{2}+it)\Gamma_q( N+\frac{3}{2}-it)  \right|\leq q^{\frac{N(N+2)}{2}} \ \ \frac{(-q^{\frac{1}{2}};q)_{\infty}}{(q^{\frac{1}{2}};q)_{\infty}}\left|  \Gamma_q(\frac{1}{2}+it)\Gamma_q( \frac{3}{2}-it)  \right|. $$
Hence,
\begin{eqnarray*}&&\left|\int_{-N-\frac{1}{2} +\frac{i\pi}{\Log(q)}}^{-N-\frac{1}{2} -\frac{i\pi}{\Log(q)} }\psi(s)ds\right|\\&\leq& q^{\frac{N(N+2)}{2}} \ \ \frac{(-q^{\frac{1}{2}};q)_{\infty}}{(q^{\frac{1}{2}};q)_{\infty}}\int_{\frac{\pi}{\Log(q)}}^{-\frac{\pi}{\Log(q)}}\left|   \Gamma_q(\frac{1}{2}+it)\Gamma_q( \frac{3}{2}-it)  \right| |\varphi(N+\frac{1}{2}-it)|x^{N+\frac{1}{2}}dt\\
&\leq& q^{\frac{N(N+2)}{2}} x^{N+\frac{1}{2}}e^{r(N+\frac{1}{2})} \frac{(-q^{\frac{1}{2}};q)_{\infty}}{(q^{\frac{1}{2}};q)_{\infty}}\int_{\frac{\pi}{\Log(q)}}^{-\frac{\pi}{\Log(q)}}\left|   \Gamma_q(\frac{1}{2}+it)\Gamma_q( \frac{3}{2}-it)  \right| \varphi_1(-t)dt.
\end{eqnarray*}
Since the function $t\mapsto \left|   \Gamma_q(\frac{1}{2}+it)\Gamma_q( \frac{3}{2}-it)  \right| \varphi_1(-t)$ is continuous on $\left[\frac{\pi}{\Log(q)},-\frac{\pi}{\Log(q)}  \right]$ and independent of $N$, and since $0<q<1$, we have
$$q^{\frac{N(N+2)}{2}} x^{N+\frac{1}{2}}e^{r(N+\frac{1}{2})} = q^{\frac{N(N+2)}{2}} q^{m(N+\frac{1}{2})}e^{r(N+\frac{1}{2})} {\rightarrow}0  \hbox{  as  } {N\rightarrow +\infty}, $$
then,
$$\lim_{N\rightarrow +\infty} \int_{-N-\frac{1}{2} -\frac{i\pi}{\Log(q)}}^{-N-\frac{1}{2} +\frac{i\pi}{\Log(q)} }\psi(s)ds=0. $$
So, taking the limit as $N\rightarrow +\infty$, in (\ref{Res}), we obtain
$$\int_{\sigma -\frac{i\pi}{\Log(q)}}^{ \sigma +\frac{i\pi}{Log(q)}}\psi(s)ds=  \frac{2i\pi (1-q)}{\Log(q)} \sum_{n=0}^\infty q^{\frac{n(n+1)}{2}}(-1)^{n}\varphi(n)x^n. $$
Finally, for all $x\in\R_{q,+}$, we have
\begin{eqnarray*} \sum_{n=0}^\infty q^{\frac{n(n+1)}{2}}(-1)^{n}\varphi(n)x^n &=&  \frac{\Log(q)}{2i\pi (1-q)}\int_{\sigma -\frac{i\pi}{\Log(q)}}^{ \sigma +\frac{i\pi}{\Log(q)}}\psi(s)ds\\ &=&  \frac{\Log(q)}{2i\pi (1-q)}\int_{\sigma -\frac{i\pi}{\Log(q)}}^{ \sigma +\frac{i\pi}{\Log(q)}}\Gamma_q(s)\Gamma_q(1-s)\varphi(-s)x^{-s}  ds.\end{eqnarray*}
Since the function $s\mapsto \Gamma_q(s)\Gamma_q(1-s)\varphi(-s) $ is analytic on $\langle 0, \delta\rangle$ and $\frac{2i\pi}{\Log(q)}$-periodic, then from Theorem \ref{car}, for all $s\in \langle 0, \delta\rangle$, we have
\begin{eqnarray*}\Gamma_q(s)\Gamma_q(1-s)\varphi(-s)&=& M_q\left[ \sum_{n=0}^\infty q^{\frac{n(n+1)}{2}}(-1)^{n}\varphi(n)x^n \right](s)\\ &=& \int_0^\infty x^{s-1}\left( \sum_{n=0}^\infty q^{\frac{n(n+1)}{2}}(-1)^{n}\varphi(n)x^n \right)d_qx. \end{eqnarray*}\endproof
{\bf Example}\\
Consider the function  defined by
$$ f(x)=E_q^{-qx}= \sum_{n=0}^\infty (-1)^{n}q^{\frac{n(n+1)}{2}}\frac{x^n}{(q;q)_n}=\sum_{n=0}^\infty (-1)^{n}q^{\frac{n(n+1)}{2}}\varphi(n)x^n,$$
with $\varphi(s)= \frac{(q^{s+1};q)_\infty}{(q;q)_\infty}$.\\
The function $\varphi$ is analytic and $\frac{2i\pi}{\Log(q)}-$ periodic on $\C$ and for $\Re(s)>-\delta>-1$, we have
$$|\varphi(s)|\leq \frac{(-1;q)_\infty}{(q;q)_\infty}.$$
Then, from the previous theorem, we have for $0<\Re(s)<\delta$,
\begin{equation}
M_q\left[E_q^{-qx}\right](s)=\Gamma_q(s)\Gamma_q(1-s)\varphi(-s)=(1-q)^s\Gamma_q(s).
\end{equation}
In particular, if the parameter $q$ satisfies $\frac{\Log(1-q)}{\Log(q)}\in\Z$, then
$$M_q\left[E_q^{-q(1-q)x}\right](s)=\Gamma_q(s).$$
A generalization of this example is given in the following result.
\begin{prop}
Let $\varphi$ be an analytic function satisfying the hypothesis of Theorem \ref{RMq2}  and $f$ be a function such that
$$f(x)=\sum_{n=0}^\infty (-1)^n q^{\frac{n(n+1)}{2}} \varphi(n) \frac{x^n}{(q ; q)_n}.$$
Then for $0<\Re(s)<\delta$,
\begin{equation}
M_q(f)(s)=(1-q)^s \varphi(-s)\Gamma_q(s).
\end{equation}
\end{prop}
\proof
It is easy to see that $$f(x)=\sum_{n=0}^\infty (-1)^n q^{\frac{n(n+1)}{2}} \varphi(n) \frac{x^n}{(q ; q)_n}=\sum_{n=0}^\infty (-1)^n q^{\frac{n(n+1)}{2}} \psi(n)x^n ,$$
with $\psi(s)= \frac{(q^{s+1};q)_\infty}{(q;q)_\infty}\varphi(s).$\\
Since $\varphi$ satisfies the conditions of Theorem \ref{RMq2}, and $s\mapsto \frac{(q^{s+1};q)_\infty}{(q;q)_\infty}$  is entire,  $\frac{2i\pi}{\Log(q)}-$ periodic on $\C$ and for $\Re(s)>-\delta>-1$, we have
$$\left|\frac{(q^{s+1};q)_\infty}{(q;q)_\infty}\right| \leq \frac{(-1;q)_\infty}{(q;q)_\infty},$$
then $\psi$ satisfies the conditions of Theorem \ref{RMq2}, and for $0<\Re(s)<1$, we have
$$M_q(f)(s)=M_q\left[ \sum_{n=0}^\infty (-1)^n q^{\frac{n(n+1)}{2}} \psi(n)x^n  \right](s)= \Gamma_q(s)\Gamma_q(1-s)\psi(-s)= (1-q)^s \varphi(-s)\Gamma_q(s).$$
\endproof
\begin{cor}
Let $\varphi$ be an analytic function on $H(\delta)$ such that the function $s \mapsto (1-q)^s \varphi(s)$ satisfies the hypothesis of Theorem \ref{RMq2} and  consider\\
$$f(x)=\sum_{n=0}^\infty (-1)^n q^{\frac{n(n+1)}{2}} \varphi(n) \frac{x^n}{[n]_q!}.$$
Then, for $0<\Re(s)<\delta$,
\begin{equation}
M_q(f)(s)=\varphi(-s)\Gamma_q(s).
\end{equation}
\end{cor}
\noindent{\bf Remark.} In the case  that $\frac{\Log(1-q)}{\Log(q)}\in \mathbb{Z}$, we have $1-q\in \R_{q,+}$. So, if a function $\varphi$   satisfies the conditions of Theorem \ref{RMq2}, then $s \mapsto (1-q)^s \varphi(s)$ satisfies them and from the previous corollary, for $0<\Re(s)<\delta$, we have
\begin{equation}
M_q\left[\sum_{n=0}^\infty (-1)^n q^{\frac{n(n+1)}{2}} \varphi(n) \frac{x^n}{[n]_q!} \right](s)=\varphi(-s)\Gamma_q(s).
\end{equation}

\begin{prop}\label{prop5}
Let $\alpha>0$, $p\in \N\cup\{0 \}$ and $\varphi$ be an analytic function on $H(\delta)$  satisfying the condition (\ref{maj})  such that the function $s\mapsto \left(1-q^\alpha\right)^s\varphi(s)$ is  $\frac{2i\pi}{\Log(q^\alpha)}-$periodic
and put
$$ f(x)= \sum_{n=0}^{\infty}(-1)^n\varphi(n)q^{\alpha \frac{n(n+1)}{2}}\frac{x^{\alpha n}}{[n+p]_{q^{\alpha}}!}.$$
Then,
\begin{equation}
M_q(f)(s)= \left[\frac{1}{\alpha}\right]_{q^{\alpha}} \varphi\left(- \frac{s}{\alpha}\right) \frac{\Gamma_{q^{\alpha}}\left(1-\frac{s}{\alpha}\right)}{\Gamma_{q^{\alpha}}\left(1-\frac{s}{\alpha}+p\right)}\Gamma_{q^{\alpha}}\left(\frac{s}{\alpha}\right).
\end{equation}
\end{prop}
\proof ~~\\
We have
$$ f(x)= \sum_{n=0}^{\infty}(-1)^n\varphi(n)q^{\alpha \frac{n(n+1)}{2}}\frac{x^{\alpha n}}{[n+p]_{q^{\alpha}}}!= h(x^{\alpha}),$$
with
$$ h(x)= \sum_{n=0}^{\infty}(-1)^n\varphi(n)q^{\alpha \frac{n(n+1)}{2}}\frac{x^{n}}{[n+p]_{q^{\alpha}}!}.$$
But, under the hypotheses of the proposition,   the function $ s\mapsto \frac{\varphi(s)}{\Gamma_{q^{\alpha}}(1+s+p)}$ is analytic, $\frac{2i\pi}{\Log(q^\alpha)}-$periodic on $H(\delta)$ and satisfies the condition (\ref{maj}). Then,  from Theorem \ref{RMq2}, we have   for $0<\Re(s)<\delta$,
$$ M_{q^{\alpha}}(h)(s)= \varphi(-s)\frac{\Gamma_{q^{\alpha}}(1-s)}{\Gamma_{q^{\alpha}}(1-s+p)}\Gamma_{q^{\alpha}}(s).$$
Finally, the relation (\ref{Mk2}) achieves the proof.\endproof
For the particular case $p=0$, if $\varphi$   satisfies  the conditions of the proposition, we get

\begin{equation}\label{CoPr5}
M_q\left[\sum_{n=0}^{\infty}(-1)^n\varphi(n)q^{\alpha \frac{n(n+1)}{2}}\frac{x^{\alpha n}}{[n]_{q^{\alpha}}!}\right]= \left[\frac{1}{\alpha}\right]_{q^{\alpha}} \varphi\left(- \frac{s}{\alpha}\right) \Gamma_{q^{\alpha}}\left(\frac{s}{\alpha}\right).
\end{equation}

\begin{prop}\label{prop6}
Let $\alpha\in \N$, $p\in \N\cup\{0 \}$ and $\varphi$ be an analytic function on $H(\delta)$  satisfying the condition (\ref{maj})  such that the function $s\mapsto \left(1-q\right)^{\alpha s}\varphi(s)$ is  $\frac{2i\pi}{\Log(q^\alpha)}-$periodic. Then, the $q$-Mellin transform of the function
$$ f(x)= \sum_{n=0}^{\infty}(-1)^n\varphi(n)q^{\alpha \frac{n(n+1)}{2}}\frac{x^{\alpha n}}{[\alpha n+p]_q!}$$
is given by
\begin{equation}
M_q(f)(s)= \left[\frac{1}{\alpha}\right]_{q^{\alpha}} \varphi\left(- \frac{s}{\alpha}\right) \frac{\Gamma_{q^{\alpha}}\left(1-\frac{s}{\alpha}\right)}{\Gamma_{q}(1-s+p)}\Gamma_{q^{\alpha}}\left(\frac{s}{\alpha}\right).
\end{equation}
\end{prop}
\proof ~~\\
We have
$$ f(x)= \sum_{n=0}^{\infty}(-1)^n\varphi(n)q^{\alpha \frac{n(n+1)}{2}}\frac{x^{\alpha n}}{[\alpha n+p]_{q}!}= h(x^{\alpha}),$$
where
$$ h(x)= \sum_{n=0}^{\infty}(-1)^n\varphi(n)q^{\alpha \frac{n(n+1)}{2}}\frac{x^{n}}{[\alpha n+p]_{q}!}.$$
Under the hypotheses of the proposition, the function $ s\mapsto \frac{\varphi(s)}{\Gamma_{q}(1+\alpha s+p)}$ is analytic, $\frac{2i\pi}{\Log(q^\alpha)}-$periodic on $H(\delta)$ and satisfies the condition (\ref{maj}). Then,  from Theorem \ref{RMq2},  we have for $0<\Re(s)<\delta$,
$$ M_{q^{\alpha}}(h)(s)= \frac{\varphi(-s)}{\Gamma_{q}(1-\alpha s+p)}\Gamma_{q^{\alpha}}(1-s)\Gamma_{q^{\alpha}}(s).$$
Hence, the result follows from the relation (\ref{Mk2}).

\section{Examples}
In this section, we shall apply the two $q$-version of the Ramanujan master theorem to find the $q$-Mellin transform of some $q$-special functions.\\
\noindent{\bf Example 1:}\\

Let $a \in \mathbb{R}_{q,+}$, $\gamma >0$ and  $f(x)= \frac{(-axq^{\gamma};q)_{\infty}}{(-ax;q)_{\infty}}.$\\

By using the $q$-Binomial formula, we get
\begin{eqnarray*}
\frac{(-axq^{\gamma};q)_{\infty}}{(-ax;q)_{\infty}}&=& \sum_{n=0}^{\infty} \frac{(q^{\gamma};q)_n}{(q;q)_n}(-ax)^n\\
&=& \sum_{n=0}^{\infty}(-1)^n \frac{\Gamma_q(\gamma+n)}{\Gamma_q(\gamma)}a^n \frac{x^n}{[n]_q!}.
\end{eqnarray*}
We consider the function $\varphi$ defined by
\begin{equation}
\varphi(s)=\frac{\Gamma_q(\gamma+s)}{\Gamma_q(\gamma)}a^s;  \qquad \Re(s)>-\gamma.
\end{equation}
Let $0<\delta< \min(1,\gamma)$. It is easy to see that the function $\varphi$ is analytic on $H(\delta)$, and the function $ s\mapsto (1-q)^s \varphi(s)$ is analytic and  $\frac{2i\pi}{\Log( q)}-$ periodic on $H(\delta)$. Furthermore, for $s=u+iv\in H(\delta)$,  we have
\begin{eqnarray*}|(1-q)^s\varphi(s)|&=&\left|(1-q)^s\frac{\Gamma_q(\gamma+s)}{\Gamma_q(\gamma)}a^s\right|=\frac{1-q}{\Gamma_q(\gamma)}|(q^{s+\gamma};q)_\infty| a^u\\ &\leq & \frac{1-q}{\Gamma_q(\gamma)}(-q^{\gamma-\delta};q)_\infty e^{u\Log(a)}.
\end{eqnarray*}
 So from Corollary \ref{cor1}, we obtain for $0<\Re(s)<\delta$,
\begin{equation}
M_q(f)(s)= \varphi(-s)\frac{\Gamma_q(s)}{K_q(s)}= \frac{\Gamma_q(\gamma-s)}{\Gamma_q(\gamma)}a^{-s}\frac{\Gamma_q(s)}{K_q(s)}=\frac{a^{-s}}{K_q(s)}B_q(\gamma-s,s).
\end{equation}
\noindent{\bf Example 2:}\\
Consider the $q$-cosine function defined by $$ \cos(x;q^2)= \sum_{n=0}^{\infty} (-1)^n q^{n(n+1)} \frac{x^{2n}}{[2n]_q!}.$$
We assume that $ \frac{\Log(1-q)}{\Log (q)}\in 2\mathbb{Z}$ and we consider the function $\varphi(s)= 1$.\\
So, the function $s\mapsto (1-q)^{2s}\varphi(s)$ is $\frac{2i\pi}{\Log (q^2)}-$ periodic.
Then by applying Proposition \ref{prop6}, for $\alpha=2$ and $p=0$, we obtain
\begin{equation}
M_q[\cos(x;q^2)](s)= \left[\frac{1}{2}\right]_{q^2} \frac{\Gamma_{q^2}\left(1-\frac{s}{2}\right)\Gamma_{q^2}\left(\frac{s}{2}\right)}{\Gamma_q(1-s)}.
\end{equation}

\noindent{\bf Example 3:}\\
We assume that $ \frac{\Log(1-q)}{\Log(q)}\in 2\mathbb{Z}$ and we consider the function $q$-sine defined by
$$ f(x)= \sin(x;q^2)= \sum_{n=0}^{\infty} (-1)^n q^{n(n+1)} \frac{x^{2n+1}}{[2n+1]_q!}=xh(x),$$
where $$h(x)=\sum_{n=0}^{\infty} (-1)^n q^{n(n+1)} \frac{x^{2n}}{[2n+1]_q!}.$$
Since, $ \frac{\Log(1-q)}{\Log(q)}\in 2\mathbb{Z}$, then  the function $s\mapsto (1-q)^{2s}$ is $\frac{2i\pi}{\Log (q^2)}-$ periodic.
So, by application of Proposition \ref{prop6}, with $\alpha=2$,   $p=1$ and $\varphi(s)=1$, we obtain
\begin{equation}
M_q(h)(s)= \left[\frac{1}{2}\right]_{q^2} \frac{\Gamma_{q^2}\left(1-\frac{s}{2}\right)\Gamma_{q^2}\left(\frac{s}{2}\right)}{\Gamma_q(2-s)}.
\end{equation}
Therefore, by using the properties of the $q$-Mellin transform, we obtain
\begin{equation}
M_q(\sin(x;q^2))(s)= M_q(h)(s+1)=\left[\frac{1}{2}\right]_{q^2} \frac{\Gamma_{q^2}\left(\frac{1-s}{2}\right)\Gamma_{q^2}\left(\frac{1+s}{2}\right)}{\Gamma_q(1-s)}.
\end{equation}

\noindent{\bf Example 4:}\\
We consider the third Jackson's $q$-Bessel function defined by
\begin{eqnarray*}
f(x)= J_{\nu}(x;q^2)&= &\frac{x^{\nu}}{(1-q^2)^{\nu} \Gamma_{q^2}(\nu+1)}~~~  _1\varphi_1(0;q^{2\nu+2};q^2;q^2x^2)\\
&=& \frac{x^{\nu}}{(1-q^2)^{\nu}}\sum_{n=0}^{\infty} (-1)^n \frac{(1-q^2)^{-2n}}{\Gamma_{q^2}(n+\nu+1)} q^{n(n+1)} \frac{x^{2n}}{[n]_{q^2}!}\\
&=& \frac{x^{\nu}}{(1-q^2)^{\nu}} h(x),
\end{eqnarray*}
where
\begin{eqnarray*}
h(x)=\sum_{n=0}^{\infty} (-1)^n \frac{(1-q^2)^{-2n}}{\Gamma_{q^2}(n+\nu+1)} q^{n(n+1)} \frac{x^{2n}}{[n]_{q^2}!}.
\end{eqnarray*}
We have
\begin{eqnarray*}
M_q(f)(s)= \frac{1}{(1-q^2)^{\nu}}M_q(x^{\nu}h(x))(s)=\frac{1}{(1-q^2)^{\nu}}M_q(h(x))(s+\nu).
\end{eqnarray*}
On the other hand,  the function $\displaystyle \varphi (s)=\frac{(1-q^2)^{-2s}}{\Gamma_{q^2}(s+\nu+1)}$,  is entire on $\C$ and for $0<\delta<1$, we have for $\Re(s)>-\delta>-1$,
 $$|\varphi (s)|=\left| (1-q^2)^{\nu-s} \frac{\left( q^{s+\nu+1};q^2 \right)_\infty}{\left(q^2;q^2  \right)_\infty} \right|\leq (1-q^2)^{\nu}\frac{\left( -q^{\nu};q^2 \right)_\infty}{\left(q^2;q^2  \right)_\infty}e^{-\Log(1-q^2)\Re(s)}.$$
 Furthermore, it is easy to verify   that the function
$ s\mapsto (1-q^2)^s \varphi (s)$ is $\frac{2i\pi}{\Log(q^2)}-$ periodic. So, by  application of the Proposition \ref{prop5} for $\alpha=2$ and $p=0$, we obtain
\begin{eqnarray}
M_q(h)(s)= \left[\frac{1}{2}\right]_{q^{2}} \varphi\left(- \frac{s}{2}\right) \Gamma_{q^{2}}\left(\frac{s}{2}\right)=\left[\frac{1}{2}\right]_{q^{2}} \frac{(1-q^2)^{s}}{\Gamma_{q^2}\left(-\frac{s}{2}+\nu+1\right)}\Gamma_{q^{2}}\left(\frac{s}{2}\right).
\end{eqnarray}
Therefore
\begin{eqnarray}
M_q\left[J_{\nu}(x;q^2)\right](s)=\left[\frac{1}{2}\right]_{q^{2}}(1-q^2)^{s}\frac{\Gamma_{q^{2}}\left(\frac{s+\nu}{2}\right)}{\Gamma_{q^2}\left(\frac{-s+\nu+2}{2}\right)}.
\end{eqnarray}
\noindent{\bf Example 5:}\\
Let
\begin{eqnarray*}
f(x)&=& _r\varphi_r (q^{a_1},q^{a_2},\dots,q^{a_r};q^{b_1},q^{b_2},\dots,q^{b_r};q;qx)\\
&=& \sum_{n=0}^{\infty}(-1)^n \frac{(q^{a_1};q)_n(q^{a_2};q)_n \dots (q^{a_r};q)_n}{(q^{b_1};q)_n(q^{b_2};,q)_n \dots (q^{b_r};q)_n }q^{\frac{n(n+1)}{2}}\frac{x^n}{(q;q)_n}\\&=& \sum_{n=0}^{\infty}(-1)^n \varphi(n)\frac{x^n}{(q;q)_n},
\end{eqnarray*}
with $$\displaystyle \varphi (s)=\frac{(q^{a_1},q)_\infty(q^{a_2},q)_\infty \dots (q^{a_r},q)_\infty}{(q^{b_1},q)_\infty(q^{b_2},q)_\infty \dots (q^{b_r},q)_\infty }\frac{(q^{b_1+s},q)_\infty(q^{b_2+s},q)_\infty \dots (q^{b_r+s},q)_\infty }{(q^{a_1+s},q)_\infty(q^{a_2+s},q)_\infty \dots (q^{a_r+s},q)_\infty}.$$
 The function  $ \varphi $ is analytic on the half-plane $\Re(s)>-\min_{1\leq k\leq r}(a_k)$  and it is $\frac{2i\pi}{\Log(q)}$ periodic. So,
\begin{equation*}\begin{split}
M_q(f)(s)&= (1-q)^s \varphi(-s) \Gamma_q(s)\\&=(1-q)^s \frac{\Gamma_q(a_1-s)\Gamma_q(a_2-s)\dots\Gamma_q(a_r-s)}{\Gamma_q(b_1-s)\Gamma_q(b_2-s)\dots\Gamma_q(b_r-s)}\frac{\Gamma_q(b_1)\Gamma_q(b_2)\dots\Gamma_q(b_r)}{\Gamma_q(a_1)\Gamma_q(a_2)\dots\Gamma_q(a_r)}\Gamma_q(s).
\end{split}\end{equation*}

\end{document}